\documentclass[]{article}
\usepackage{amssymb,amsfonts,amsmath}%
\usepackage{latexsym}
\usepackage{graphicx,epsfig}

\usepackage{epsfig}

\newtheorem{thm}{Theorem}
\newtheorem{pro}[thm]{Proposition}
\newtheorem{lem}[thm]{Lemma}
\newtheorem{cl}[thm]{Claim}

\newcommand{\prf}{\noindent{\it Proof.} }
\newcommand {\cbdo}{\hfill$\Box$}

\begin{document}

\title{An asymptotic solution to the Erd\H{o}s four-edge intersection problem}
\author{Andrzej \.{Z}ak\thanks{The author was partially supported
by the Polish Ministry of Science and Higher Education.}\\
\small{AGH University of Krakow,
Poland}} \maketitle

\begin{abstract}
For an $n$-vertex graph $G$ and a permutation $\sigma$ of its vertex set,
let $\sigma(G)$ denote the corresponding relabelling of $G$, and put
$I_G(\sigma)=|E(G)\cap E(\sigma(G))|$. Let $f(n,k)$ be the
minimum number of edges in an $n$-vertex graph for which
$I_G(\sigma)\geq k$ for every $\sigma$.
In his 1977 formulation of the problem, Erd\H{o}s discussed the small
values of $k$ and left the cases $k=4$ and $k=5$ as the next natural
open questions.  For $k=4$ he asked whether $f(n,4)=2n-4$, with the
upper bound witnessed by $K_{2,n-2}$; the neighbouring $k=5$ question
was recently settled exactly by Fang and Hou.
We prove that every graph $G$ of order $n$ and size at most
$2n-10n^{2/3}-7$ has a relabelling with at most three common edges.
Consequently,
\[
 2n-10n^{2/3}-7<f(n,4)\leq 2n-4,
\]
and hence
\[
 f(n,4)=2n-o(n).
\]
Thus we resolve Erd\H{o}s's four-edge intersection problem asymptotically,
confirming his proposed value up to a sublinear error term.  For comparison, for all sufficiently large $n$, Fang and Hou's result
guarantees at most four common edges for graphs with at most $2n-3$ edges,
whereas reducing the edge bound by only $10n^{2/3}+4=o(n)$ already allows
us to guarantee at most three common edges.
\end{abstract}

\section{Introduction}
\begingroup

In this paper, a graph is simple and has no loops or multiple edges.  For
a graph $G$, its vertex and edge sets are denoted by $V(G)$ and $E(G)$,
respectively.  If $\sigma$ is a permutation of $V(G)$, then $\sigma(G)$
denotes the graph with vertex set $V(G)$ and edge set
\[
 E(\sigma(G))=\{\sigma(u)\sigma(v):uv\in E(G)\}.
\]
A permutation $\sigma$ is a packing of $G$ if $G$ and $\sigma(G)$ are
edge-disjoint.  More generally, put
\[
 I_G(\sigma):=|E(G)\cap E(\sigma(G))|.
\]
Equivalently,
\[
 I_G(\sigma)=
 \bigl|\{uv\in E(G):\sigma(u)\sigma(v)\in E(G)\}\bigr|.
\]

In 1977, Erd\H{o}s~\cite{Erdos} introduced the following minimum-intersection
problem.  For $k\geq1$, let
\[
 f(n,k):=\min\bigl\{|E(G)|:|V(G)|=n\text{ and }
 I_G(\sigma)\geq k\text{ for every permutation }\sigma\bigr\}.
\]
Thus every $n$-vertex graph with fewer than $f(n,k)$ edges has a
relabelling sharing at most $k-1$ edges with the original graph.

Erd\H{o}s then considered the small values of $k$.  He recorded
$f(n,1)=n-1$ and, citing Chung, Graham, and Murty,
\[
 f(n,2)=f(n,3)=\left\lfloor\frac{3n}{2}\right\rfloor,
\]
and left the next two cases, $k=4$ and $k=5$, open.  For $k=4$ he observed
\begin{equation}\label{upper-four}
 f(n,4)\leq 2n-4\qquad(n\geq6),
\end{equation}
using the complete bipartite graph $K_{2,n-2}$, and asked whether
\begin{equation}\label{erdos-conjecture}
 f(n,4)=2n-4.
\end{equation}
Immediately afterwards he recorded Mullin's construction giving
$f(n,5)\leq2n-2$ and suggested equality there as well.  Thus the four-
and five-edge questions appeared side by side as the first unresolved small
cases in Erd\H{o}s's formulation.  We refer to~\eqref{erdos-conjecture} as
the Erd\H{o}s four-edge intersection problem and determine its asymptotic
answer.

The minimum-intersection viewpoint was subsequently placed in a broader
framework by Mullin, Roy, and Schellenberg~\cite{MRS}.  We shall use the
closely related language of near-packings.  Given a family of graphs
$\mathcal F$, a permutation $\sigma$ is called an $\mathcal F$-near-packing
of $G$ if the graph formed by the common edges
$E(G)\cap E(\sigma(G))$ belongs to $\mathcal F$.  Let
$\mu(n,\mathcal F)$ be the maximum integer $m$ such that every graph of
order $n$ and size at most $m$ has an $\mathcal F$-near-packing.  This
notion was introduced in~\cite{zak1}; an earlier bounded-degree version
appeared in~\cite{E}, and near-packings of two not necessarily
isomorphic graphs were considered in~\cite{KZ}.

For $r\geq0$, let $\mathcal E_r$ denote the family of all graphs with at
most $r$ edges.  Directly from the definitions,
\begin{equation}\label{relation-f-mu}
 f(n,k)=\mu(n,\mathcal E_{k-1})+1.
\end{equation}
Hence Erd\H{o}s's four-edge problem is precisely the problem of determining
$\mu(n,\mathcal E_3)$.

The neighbouring case $k=5$ has very recently been settled exactly by
Fang and Hou~\cite{FH}, who proved $f(n,5)=2n-2$ for every sufficiently
large $n$.  Equivalently, every such graph with at most $2n-3$ edges has
a relabelling with at most four common edges.  Our theorem lowers the edge
bound by only $10n^{2/3}+4=o(n)$, but already improves the guaranteed
intersection from four common edges to three.  Thus the two results operate
at essentially the same linear edge scale $2n$, while differing by one in
the guaranteed intersection.

Our main result is the following.
\begin{thm}\label{E3main}
For every $n\geq6$,
\[
 2n-10n^{2/3}-7<f(n,4)\leq2n-4.
\]
Consequently,
\[
 f(n,4)=2n-o(n).
\]
Equivalently,
\[
 \mu(n,\mathcal E_3)\geq
 \left\lfloor2n-10n^{2/3}-7\right\rfloor.
\]
\end{thm}
The upper bound in Theorem~\ref{E3main} is Erd\H{o}s's construction
$K_{2,n-2}$.  The remainder of the paper is devoted to the lower bound,
which proves~\eqref{erdos-conjecture} asymptotically.
\endgroup


\section{Preliminaries}
\begingroup

\begin{thm}[\cite{BE,BS,SS}]\label{pak1}
Every graph of order $n$ and size at most $n-2$ is packable.
\end{thm}
\endgroup

The following lemma has become a standard technique for extending a packing that covers nearly the entire graph to a packing of the entire graph. A proof based on Hall's Theorem can be found, for instance, in \cite{GPWZ, GZ2a}.
\begin{lem}\label{duzomalych}
Let $G$ be a graph and $k \geq 1$ be any integer. Suppose that
$G$ contains an independent set $U$ such that
\begin{enumerate}
\item $d_G(u)\leq k$ for each $u\in U$,
\item $N_G(u)\cap N_G(v) = \emptyset$ for every distinct $u,v \in U$.
\end{enumerate}
If $|U|\geq 2k$, then for every permutation $\sigma'$ of $V(G)\setminus U$
there exists a permutation $\sigma$ of $V(G)$ such that
\[
E(\sigma'(G-U))\cap E(G-U) = E(\sigma(G))\cap E(G).
\]
\end{lem}
\prf Let $G':=G-U$ and $\sigma'$ be any permutation of $V(G')$. Below we
show that we can extend $\sigma'$ to a permutation $\sigma$ as required of $G$.\\
For any $v \in V(G')$ let us define $\sigma(v):=\sigma'(v)$. Then
let us consider a bipartite graph $B$ with partition sets
$X:=U \times \{0\}$ and $Y:=U \times
\{1\}$. For $u,v \in U$ the vertices $(u,0)$,
$(v,1)$ are joined by an edge in $B$ if and only if
$\sigma'(N_G(u)) \cap N_G(v) = \emptyset$. So, if $(u,0)$, $(v,1)$
are joined by an edge in $B$ we can put $\sigma(u)=v$.
Therefore, since $d_G(u) \leq k$ for $u \in U$, we have $d_B((u,0))
\geq |U|-k\geq |U|/2$, by the assumption on $|U|$. Similarly,
$d_B((v,1)) \geq |U|/2$.

Let $S \subset X$. If $|S| \leq |U|/2$ then obviously $|N_B(S)| \geq |S|$. Notice
that if $|S| >|U|/2$ then $N_B(S)=Y$. Indeed, otherwise let $(v,1)\in Y$ be a vertex which has
no neighbour in $S$. Thus,
\begin{align*}
 d_B((v,1))\leq {|U|}-|S|< |U|- |U|/2= |U|/2,
\end{align*}
a contradiction. Hence,
in any case $|S|\leq |N_B(S)|$. Thus, by Hall's
theorem there is a perfect matching $M$ in $B$. Therefore we can define
$\sigma(u)=v$ for $u,v\in U$ such that $(u,0)$,
$(v,1)$ are incident with the same edge in $M$. \cbdo

\begingroup

\begin{pro}\label{indukcja}
Let $G$ be a graph of order $n$ and size $m$ with $m\leq an-g(n)$, where $a$ is a real number
and $g(n)$ is a non-decreasing function.
If $U\subset V(G)$ and vertices from $U$ cover at least $a|U|$ edges, then
\[m'\leq an'-g(n'),\]
where $n'$ and $m'$ are respectively the order and the size of $G-U$.
\end{pro}
\prf
\begin{align*}
m'&\leq an-g(n)-a|U|=a(n-|U|)-g(n)\\
&\leq a(n-|U|)-g(n-|U|)=an'-g(n'),
\end{align*}
because $g(n)\geq g(n-|U|)$.
\cbdo

\begin{thm}[\cite{SS}]\label{Sp}
Let $G_1$ and $G_2$ be graphs of order $n$.  If
\[
2\Delta(G_1)\Delta(G_2)<n,
\]
then the complete graph $K_n$ contains edge-disjoint copies of $G_1$ and
$G_2$.
\end{thm}
\endgroup

\begingroup

\begin{lem}\label{globalE3}
Let $G$ be a graph, let $s\geq1$, and let
\[
C=\{x_1,\ldots,x_s\}\subseteq V(G),
\]
and let $T_1,\ldots,T_{2s}$ be distinct components of $G-C$.
Suppose that every $T_j$ is a tree and
\[
d_G(v,C)\leq 1
\]
for every
\[
v\in V(T_1)\cup\cdots\cup V(T_{2s}).
\]
Put
\[
W:=C\cup V(T_1)\cup\cdots\cup V(T_{2s})
\qquad\text{and}\qquad
H:=G[W].
\]
Then there exists a permutation $\sigma$ of $W$ such that
\[
|E(H)\cap E(\sigma(H))|\leq 3
\]
and
\[
\sigma(C)\subseteq W\setminus C,
\qquad
\sigma^{-1}(C)\subseteq W\setminus C.
\]
In particular, $H$ has an $\mathcal{E}_3$-near-packing which
moves every vertex of $C$ outside $C$ in both directions.
\end{lem}

\prf
Pair the tree components arbitrarily, and for $i\in[s]$ put
\[
P_i:=G[V(T_{2i-1})\cup V(T_{2i})].
\]
Thus $P_i$ is the disjoint union of two trees.

For every $i\in[s]$, choose a vertex $a_i\in V(P_i)$ as follows.
If at least one of the two trees in $P_i$ is non-trivial, let
$a_i$ be a leaf of a non-trivial tree. If both trees are trivial,
choose either of their vertices as $a_i$. Put
\[
F_i:=P_i-a_i.
\]
If $a_i$ belongs to a non-trivial tree, denote by $b_i$ its unique
neighbour in $P_i$.

Unless both trees in $P_i$ are trivial, the graph $F_i$ is the
disjoint union of two non-empty trees. Consequently,
\[
|E(F_i)|=|V(F_i)|-2.
\]
Hence, by Theorem~\ref{pak1}, there exists a proper self-packing
$\pi_i$ of $F_i$. If both trees in $P_i$ are trivial, then $F_i$
consists of one vertex, and we take $\pi_i$ to be the identity
permutation.

For every vertex
\[
v\in V(T_1)\cup\cdots\cup V(T_{2s}),
\]
define its \emph{core label} by
\[
\ell(v):=
\begin{cases}
r, & \text{if } vx_r\in E(G),\\
\bot, & \text{if } N_G(v)\cap C=\emptyset.
\end{cases}
\]
This is well-defined because $d_G(v,C)\leq 1$.

We first choose a permutation $\alpha$ of $[s]$. Whenever $b_i$
is defined, put
\[
w_i:=\pi_i^{-1}(b_i)
\qquad\text{and}\qquad
\lambda_i:=\ell(w_i).
\]
For a uniformly random permutation $\alpha\in S_s$, let
\[
X(\alpha):=
\sum_{\substack{i\in[s]\\ b_i\text{ is defined}\\
                         \lambda_i\neq\bot}}
\mathbf{1}_{\{\alpha(\lambda_i)=i\}}.
\]
For every term occurring in this sum,
\[
\Pr\bigl(\alpha(\lambda_i)=i\bigr)=\frac1s.
\]
There are at most $s$ terms, and hence
\[
\mathbb{E}X(\alpha)\leq 1.
\]
We may therefore fix a permutation $\alpha\in S_s$ such that
\begin{equation}\label{alpha-conflicts}
X(\alpha)\leq 1.
\end{equation}

We next choose a permutation $\beta$ of $[s]$. For every $i$ for
which $b_i$ is defined, put
\[
\rho_i:=\ell(\pi_i(b_i)).
\]
Also put
\[
\delta_i:=\ell(a_i).
\]
After fixing $\alpha$, define
\[
\theta_i:=
\begin{cases}
\delta_{\alpha(\delta_i)},
   & \text{if }\delta_i\neq\bot
     \text{ and }\delta_{\alpha(\delta_i)}\neq\bot,\\
\bot, & \text{otherwise}.
\end{cases}
\]
For a uniformly random permutation $\beta\in S_s$, let
\[
Y(\beta):=
\sum_{\substack{i\in[s]\\ b_i\text{ is defined}\\
                         \rho_i\neq\bot}}
\mathbf{1}_{\{\beta(i)=\rho_i\}}
+
\sum_{\substack{i\in[s]\\ \theta_i\neq\bot}}
\mathbf{1}_{\{\beta(i)=\theta_i\}}.
\]
Each indicator has expectation $1/s$, and there are at most $2s$
of them. Therefore
\[
\mathbb{E}Y(\beta)\leq 2.
\]
We may thus fix a permutation $\beta\in S_s$ such that
\begin{equation}\label{beta-conflicts}
Y(\beta)\leq 2.
\end{equation}

We now define a permutation $\sigma$ of $W$. For $r,i\in[s]$, let
\[
\sigma(x_r):=a_{\alpha(r)},
\qquad
\sigma(a_i):=x_{\beta(i)},
\]
and, for every $v\in V(F_i)$, let
\[
\sigma(v):=\pi_i(v).
\]
This is indeed a permutation: the three sets
\[
C,\qquad A:=\{a_1,\ldots,a_s\},
\qquad \bigcup_{i=1}^{s}V(F_i)
\]
are mapped bijectively onto, respectively,
\[
A,\qquad C,\qquad \bigcup_{i=1}^{s}V(F_i).
\]
In particular,
\[
\sigma(C)=A=\sigma^{-1}(C),
\]
so both $\sigma(C)$ and $\sigma^{-1}(C)$ are disjoint from $C$.

It remains to count the common edges. First, no edge with both
endvertices in some $F_i$ becomes a common edge, since $\pi_i$ is
a proper self-packing of $F_i$. Moreover, there are no edges
between different $P_i$'s. An edge with both endvertices in $C$
is mapped to a pair
\[
a_{\alpha(r)}a_{\alpha(r')}.
\]
The two vertices belong to different components of $G-C$, so this
pair is not an edge.

Consider an edge $x_rv$, where $v\in V(F_i)$. Its image is
\[
a_{\alpha(r)}\pi_i(v).
\]
If $\alpha(r)\neq i$, its endpoints belong to different components
of $G-C$, so it is not an edge. If $\alpha(r)=i$, then it can be an
edge only when $b_i$ is defined and
\[
\pi_i(v)=b_i.
\]
In that case $v=w_i$, and the original edge $x_rv$ implies
\[
r=\ell(w_i)=\lambda_i.
\]
Thus all common edges of this type are counted by $X(\alpha)$.
By~\eqref{alpha-conflicts}, there is at most one such edge.

The only edge of $P_i$ incident with $a_i$ is $a_ib_i$, when
$b_i$ is defined. Its image is
\[
x_{\beta(i)}\pi_i(b_i).
\]
This is an edge precisely when
\[
\beta(i)=\ell(\pi_i(b_i))=\rho_i.
\]
Such common edges are counted by the first sum in $Y(\beta)$.

Finally, suppose that $a_i$ is adjacent to a vertex of $C$.
Then this vertex is $x_{\delta_i}$, and
\[
\sigma(x_{\delta_i}a_i)
 =
a_{\alpha(\delta_i)}x_{\beta(i)}.
\]
This is an edge exactly when
\[
\beta(i)=\ell(a_{\alpha(\delta_i)})
        =\delta_{\alpha(\delta_i)}
        =\theta_i.
\]
These common edges are counted by the second sum in $Y(\beta)$.

The above cases exhaust all edges of $H$. Consequently,
by~\eqref{alpha-conflicts} and~\eqref{beta-conflicts},
\[
|E(H)\cap E(\sigma(H))|
 \leq X(\alpha)+Y(\beta)
 \leq 1+2=3.
\]
Thus $\sigma$ is the required $\mathcal{E}_3$-near-packing of
$H$.
\cbdo
\endgroup

\section{Proof of Theorem \ref{E3main}}
\begingroup

Let
\[
t=\left\lfloor n^{1/3}\right\rfloor
\qquad\text{and}\qquad
g(n)=10n^{2/3}+7.
\]
We prove that every graph of order $n$ and size at most
\[
2n-g(n)=2n-10n^{2/3}-7
\]
has an $\mathcal E_3$-near-packing.  Suppose that $G$ is a
counterexample of minimum order $n$.

If $n\leq1000$, then $10n^{2/3}\geq n$, and hence
\[
|E(G)|\leq n-7\leq n-2.
\]
Theorem~\ref{pak1} would give a proper packing of $G$.  Consequently,
$n>1000$.  By Theorem~\ref{Sp},
\[
2\Delta(G)^2\geq n.
\]
Thus $\Delta(G)\geq23$.  Fix a vertex $u$ with
$d_G(u)=\Delta(G)=:\Delta$.

For a set $D\subseteq V(G)$, let $e_G(D)$ denote the number of edges
having at least one endpoint in $D$.  Whenever $e_G(D)\geq2|D|$,
Proposition~\ref{indukcja}, with $a=2$, and the minimality of $G$ imply
that $G-D$ has an $\mathcal E_3$-near-packing.

\begin{cl}\label{few-leaves}
The graph $G$ has no isolated vertices and has at most seven vertices
of degree one.
\end{cl}
\prf
Suppose first that $v$ is isolated.  The set $\{u,v\}$ covers
$\Delta\geq4$ edges, so $G-\{u,v\}$ has an $\mathcal E_3$-near-packing
$\tau$.  Extend $\tau$ by the transposition $(u\ v)$.  Every edge
incident with $u$ is mapped to a pair incident with the isolated vertex
$v$, so no new common edge is created.  This contradicts the choice of
$G$.

Assume now that $G$ has at least eight vertices of degree one, and let
$L$ be the set of all such vertices.  If two vertices $v,w\in L$ are
adjacent, then $vw$ is a component of $G$.  The set $\{u,v,w\}$ covers
$\Delta+1\geq6$ edges.  An $\mathcal E_3$-near-packing of
$G-\{u,v,w\}$ extends, without creating a common edge meeting this set,
by the cycle
\[
v\mapsto w\mapsto u\mapsto v.
\]
Hence we may assume that $L$ is independent.

For $v\in L$, denote its unique neighbour by $y(v)$, and put
$Y:=\{y(v):v\in L\}$.  We distinguish three cases.

\medskip
\noindent\emph{Case 1. $|Y|=1$.}
Let $y$ be the common neighbour of all vertices of $L$, and choose
distinct $p,q\in L$.  The set $D=\{y,p,q\}$ covers $d(y)\geq8>2|D|$
edges.  Put $R=G-D$, and let $\tau$ be an $\mathcal E_3$-near-packing
of $R$.  If $\tau$ has three common edges, choose an edge $ab\in E(R)$
such that
\[
\tau(a)\tau(b)\in E(R),
\]
and fix one of its endpoints, say $a$; otherwise choose $a\in V(R)$
arbitrarily.  Put $a^*:=\tau(a)$ and define
\[
\sigma(y)=p,\qquad \sigma(a)=q,\qquad
\sigma(p)=a^*,\qquad \sigma(q)=y,
\]
while $\sigma(x)=\tau(x)$ for $x\in V(R)\setminus\{a\}$.
This is a permutation of $V(G)$.  Every edge of $R$ incident with $a$
is mapped to a pair incident with $q$, and hence to a non-edge, since
the only neighbour of $q$ is $y\notin V(R)$.  Thus, if $\tau$ has
three common edges, at least one of the three old overlaps disappears.
Every edge of $R$ not incident with $a$ retains its previous image.  The
only possible new common edge is the image of $yq$, namely $py$.
Consequently, $\sigma$ has at most three common edges, a contradiction.

\medskip
\noindent\emph{Case 2. $|Y|\geq2$ and some $y\in Y$ is adjacent to at
least two leaves.}
Choose distinct leaves $p,q$ adjacent to $y$, and choose a leaf $w$
with a different neighbour $z$.  If $u\in\{y,z\}$, then
$D=\{p,y,w,z\}$ covers at least $\Delta$ edges and hence at least
$2|D|$ edges.  An $\mathcal E_3$-near-packing of $G-D$ extends without
new common edges by the cycle
\[
(p\ y\ w\ z).
\]
If $u\notin\{y,z\}$, put
\[
D=\{u,y,z,p,q,w\}.
\]
The vertices $u,y,z$ cover at least
\[
d(u)+d(y)+d(z)-3\geq\Delta\geq2|D|
\]
edges.  An $\mathcal E_3$-near-packing of $G-D$ extends without new
common edges by the cycle
\[
(u\ p\ y\ w\ z\ q).
\]
Indeed, each of $u,y,z$ is mapped to a leaf whose only neighbour lies
in $D$, and all edges inside $D$ are mapped to non-edges.  Both
possibilities give a contradiction.

\medskip
\noindent\emph{Case 3. Distinct leaves have distinct neighbours.}
Choose eight leaves and denote their distinct neighbours by
$y_1,\ldots,y_8$.  If for some $i\neq j$,
\[
d(y_i)+d(y_j)-\mathbf 1_{\{y_iy_j\in E(G)\}}\geq8,
\]
then the four vertices consisting of $y_i,y_j$ and their leaves cover
at least eight edges.  As in Case~2, the corresponding four-cycle
extends an $\mathcal E_3$-near-packing of the remaining graph without
creating a new common edge.

We may therefore assume that no such pair exists.  Then some $y_i$ has
degree at most three.  Indeed, otherwise every non-adjacent pair among
the $y_i$'s would have degree sum at least eight, so the eight vertices
$y_1,\ldots,y_8$ would form a clique.  But then every $y_i$ would have
degree at least eight, which is again impossible under the displayed
inequality.

Let $z$ be a neighbour of degree at most three and let $w$ be its leaf.
Among the other seven leaves we can choose five, say
$v_1,\ldots,v_5$, such that
\[
N(v_i)\cap N(z)=\emptyset\qquad(i=1,\ldots,5),
\]
because at most $d(z)-1\leq2$ of their distinct neighbours belong to
$N(z)$.  Put
\[
U:=\{z,v_1,\ldots,v_5\}.
\]
The set $U$ is independent, its vertices have degree at most three,
and their neighbourhoods are pairwise disjoint.  Moreover, $|U|=6$.
The set
\[
D:=U\cup\{u,w\}
\]
has eight vertices and covers at least
\[
d(u)+d(z)+5-2\geq\Delta+4\geq16=2|D|
\]
edges.  Hence $G-D$ has an $\mathcal E_3$-near-packing.  Extend it to
$G-U$ by transposing $u$ and $w$.  Since $z\in U$, the vertex $w$ is
isolated in $G-U$, and this transposition creates no new common edge.
Finally, Lemma~\ref{duzomalych}, applied with $k=3$, extends this
permutation over $U$ without changing the set of common edges.  This
last contradiction proves the claim. \cbdo

By Claim~\ref{few-leaves}, every vertex of $G$ has degree at least two,
with at most seven exceptions of degree one.  Let $q$ be the number of
vertices of degree greater than $t$.  Since every such degree is at
least $n^{1/3}$,
\begin{align*}
4n-20n^{2/3}-14
 &\geq 2|E(G)|=\sum_{v\in V(G)}d_G(v)\\
 &\geq 2n-7+q\bigl(n^{1/3}-2\bigr).
\end{align*}
Consequently, either the numerator below is negative, which is already a
contradiction, or
\[
q\leq\frac{2n-20n^{2/3}-7}{n^{1/3}-2}<2n^{2/3}.
\]

Choose a maximal set $S\subseteq V(G)$ such that
\begin{enumerate}
\item $S$ is independent;
\item $2\leq d_G(s)\leq t$ for every $s\in S$;
\item the vertices of $S$ have pairwise disjoint neighbourhoods.
\end{enumerate}
The preceding bound implies $S\neq\emptyset$.  Since the vertices of
$S$ cover at least $2|S|$ edges, $G-S$ has an
$\mathcal E_3$-near-packing.  If $|S|\geq2t$, Lemma~\ref{duzomalych}
would extend this near-packing to $G$.  Hence
\[
|S|<2t.
\]
Put
\[
C:=N_G(S),\qquad c:=|C|.
\]
Then
\begin{equation}\label{C-bound}
c<2t^2\leq2n^{2/3}.
\end{equation}

For $j=2,\ldots,t$, let
\[
V_j:=\{v\in V(G)\setminus C:d_G(v)=j\}.
\]
By the maximality of $S$, every vertex in $V_2\cup\cdots\cup V_t$ has
a neighbour in $C$.  For each $v\in V_2\cup\cdots\cup V_t$, choose
one edge $e_v$ joining $v$ to $C$.  These chosen edges are pairwise
distinct, and hence their number is at least
\begin{align}
|V_2\cup\cdots\cup V_t|
 &\geq n-7-q-c
 >n-7-4n^{2/3}.                                      \label{coverns}
\end{align}

Consider the graph $G-C$.  Let $T_1,\ldots,T_p$ be those components of
$G-C$ which are trees and in which every vertex has at most one
neighbour in $C$.  We call them the \emph{minimal components} of
$G-C$.  Let $R$ be the subgraph induced by all remaining components,
and put
\[
r:=|E(R)|+\bigl|\{v\in V(R):d_G(v,C)\geq2\}\bigr|.
\]
Clearly $r\geq|E(R)|$.  Also $r\geq|V(R)|$: a non-tree component has at
least as many edges as vertices, whereas every tree component of $R$
contains a vertex with at least two neighbours in $C$.

There are exactly
\[
n-c-|V(R)|-p
\]
edges in $T_1\cup\cdots\cup T_p$.  Moreover, all edges of $R$ are
disjoint from the chosen edges $e_v$.  For every
$v\in V(R)$ with $d_G(v,C)\geq2$, there is also an edge from $v$ to
$C$ different from $e_v$ whenever $e_v$ was chosen.  Thus, in addition
to the edges counted in~\eqref{coverns}, the quantity $r$ accounts for
all edges of $R$ and for at least one further edge from $C$ for every
vertex of $R$ having at least two neighbours in $C$.  Hence
\begin{align*}
2n-10n^{2/3}-7
 &\geq |E(G)|\\
 &\geq n-7-4n^{2/3}
       +\bigl(n-c-|V(R)|-p\bigr)+r\\
 &>2n-6n^{2/3}-7-p-|V(R)|+r.
\end{align*}
It follows that
\begin{equation}\label{p}
p>4n^{2/3}-|V(R)|+r.
\end{equation}
Since $c<2t^2$ and $c$ is an integer,
\[
2c+2\leq4t^2\leq4n^{2/3}.
\]
Combining this with~\eqref{p} gives
\begin{equation}\label{p-strong}
p>2c+2-|V(R)|+r\geq2c+2.
\end{equation}

Let
\[
G':=G\bigl[C\cup V(T_1)\cup\cdots\cup V(T_{2c})\bigr]
\qquad\text{and}\qquad
G'':=G-G'.
\]
By Lemma~\ref{globalE3}, $G'$ has an $\mathcal E_3$-near-packing
$\sigma'$ such that
\[
\sigma'(C)\cap C=\emptyset
\qquad\text{and}\qquad
(\sigma')^{-1}(C)\cap C=\emptyset.
\]

Using~\eqref{p-strong},
\begin{align*}
|E(G'')|
 &=|E(R)|+\sum_{i=2c+1}^{p}\bigl(|V(T_i)|-1\bigr)\\
 &<|E(R)|+\sum_{i=2c+1}^{p}|V(T_i)|
       -(r-|V(R)|)-2\\
 &\leq |V(R)|+\sum_{i=2c+1}^{p}|V(T_i)|-2\\
 &=|V(G'')|-2.
\end{align*}
Thus, by Theorem~\ref{pak1}, $G''$ has a proper packing $\sigma''$.

Define $\sigma$ on $V(G)$ by using $\sigma'$ on $V(G')$ and
$\sigma''$ on $V(G'')$.  There are at most three common edges inside
$G'$ and none inside $G''$.  It remains to exclude a common edge joining
the two parts.  Every edge from $G'$ to $G''$ has its endpoint in $G'$
inside $C$, because the $T_i$ are components of $G-C$.  If such an edge
were mapped to another edge between the two parts, its endpoint in
$C$ would also have its image in $C$, contradicting
$\sigma'(C)\cap C=\emptyset$.  Therefore $\sigma$ is an
$\mathcal E_3$-near-packing of $G$, contrary to the choice of $G$.

No counterexample exists.  Hence every graph with
$|E(G)|\leq2n-10n^{2/3}-7$ has a relabelling with at most three common
edges.  In view of~\eqref{relation-f-mu} and~\eqref{upper-four},
Theorem~\ref{E3main} follows. \cbdo
\endgroup

\end{document}